\documentclass[12pt]{article}
\usepackage{amsmath,amssymb,amsbsy,amsfonts,amsthm,latexsym,
                 amsopn,amstext,amsxtra,euscript,amscd}

\def\mand{\qquad\mbox{and}\qquad}

\def\\{\cr}
\def\({\left(}
\def\){\right)}
\def\[{\left[}
\def\]{\right]}
\def\<{\langle}
\def\>{\rangle}
\def\fl#1{\left\lfloor#1\right\rfloor}

\begin{document}

\title{On Some Weighted Average Values of $L$-functions}

\author{ {\sc Igor E.~Shparlinski\thanks{During the preparation of this work,
the author was supported in part by 
ARC Grant DP0556431}} \\
{Department of Computing, 
Macquarie University} \\
{Sydney, NSW 2109, Australia} \\
{igor@ics.mq.edu.au}}

\date{}

\maketitle

\begin{abstract}
Let $q\ge 2$ and $N\ge 1$ be integers.
 W.~Zhang (2008) has shown that for any fixed $\varepsilon> 0$, 
 and    $q^{\varepsilon} \le N \le q^{1/2 -\varepsilon}$,  
 $$
  \sum_{\chi \ne \chi_0} \left|\sum_{n=1}^N \chi(n)\right|^2
|L(1, \chi)|^2 = (1  + o(1)) \alpha_q q N, 
$$
where the sum is take over all nonprincipal characters $\chi$ modulo  $q$,
$L(s, \chi)$ is the  $L$-functions $L(1, \chi)$ corresponding to $\chi$ 
and  $\alpha_q = q^{o(1)}$ is some explicit function of $q$. 
Here we improve this result and show that the same asymptotic 
formula holds in the essentially full range 
$q^{\varepsilon} \le N \le q^{1 -\varepsilon}$.
   \end{abstract}



\newpage 

\section{Introduction}

For integers $q\ge 2$ and $N\ge 1$ we consider the average value
$$
S(q;N) = \sum_{\chi \ne \chi_0} \left|\sum_{n=1}^N \chi(n)\right|^2
|L(1, \chi)|^2
$$
taken over all nonprincipal characters $\chi$ modulo an integer $q \ge 2$,
with $L$-functions $L(1, \chi)$ corresponding to $\chi$, weighted 
by incomplete character sums.

W.~Zhang~\cite{Zhang} has given 
an asymptotic formula for $S(q;N)$ that 
is nontrivial for 
$q^{\varepsilon} \le N \le q^{1/2 -\varepsilon}$ for any 
fixed $\varepsilon > 0$ and sufficiently large $q$. 
 
Here we improve the error term of that formula which makes it 
nontrivial in the range 
$q^{\varepsilon} \le N \le q^{1 -\varepsilon}$.

More precisely, let 
$$
\alpha_q = \(\beta_q + \gamma_q\) \frac{\varphi(q)^2}{q^2},
$$
where
\begin{eqnarray*}
\beta_q & = & \frac{\pi^2}{6}  
 \prod_{\substack{p\mid q\\p~\textrm{prime}}}\(1 - \frac{1}{p^2}\),\\
 \gamma_q & = & \frac{\pi^2}{3 \zeta(3)}  
 \prod_{\substack{p\mid q\\p~\textrm{prime}}}\(1 - \frac{1}{p^2+p+1}\)
 \sum_{\substack{m,n=1\\\gcd(nm(n+m),q)=1}}\frac{1}{nm(n+m)},
\end{eqnarray*}
$\zeta(s)$ is the Riemann zeta-function 
and $\varphi(q)$ denotes the Euler function.

It is shown in~\cite{Zhang} that
 \begin{equation}
 \label{eq:Zhang}
S(q,N) = \alpha_q q N + 
O\(\varphi(q) 2^{\omega(q)} (\log q)^2 + N^3 (\log q)^2 \),
\end{equation}
where $\omega(q)$ is the number of prime divisors of $q$. 

Since 
 \begin{equation}
 \label{eq:phi tau}
2^{\omega(q)} \le \tau(q) = q^{ o(1)}  
 \mand \varphi(q) = q^{1 + o(1)} ,
\end{equation}
where $\tau(q)$ is the number of positive integer divisors of $q$, 
see~\cite[Theorems~317 and~328]{HardyWright},
we conclude that $ \alpha_q = q^{ o(1)}$ and the error in~\eqref{eq:Zhang}
is of the shape $O\( q^{1 + o(1)} + N^3  q^{o(1)}\)$.

In particular, the asymptotic formula~\eqref{eq:Zhang} 
 is nontrivial if
 $q^{\varepsilon} \le N \le q^{1/2 -\varepsilon}$ for any 
 fixed $\varepsilon > 0$ and $q$ is  large  enough.

 Here we estimate a certain sum which arises in~\cite{Zhang}
 in a more accurate way and 
 essentially replace $N^3$ in~\eqref{eq:Zhang}
 with $N^2q^{o(1)}$ which makes it nontrivial in the range 
   $q^{\varepsilon} \le N \le q^{1 -\varepsilon}$.

\section{Main Result}

{\bf Theorem.}\ {\it 
Let $q > N \ge 1$ be integers. Then  
$$
S(q,N) = \alpha_q q N + 
O\(\varphi(q) 2^{\omega(q)} (\log q)^2 + N^2q^{o(1)} \), 
$$
as $q \to \infty$. 
\/}

\begin{proof} It has been shown in~\cite{Zhang}
that 
$$
S(q,N) = M_1 + M_2 + O\(N^2 (\log q)^2\), 
$$
where 
$$
M_1 =\varphi(q)\sum_{m,n=1}^N \sum_{\substack{u,v=1\\ mu = nv }}^{q^2}\frac{1}{uv}
\mand
M_2 =\varphi(q)\sum_{m,n=1}^N \sum_{\substack{u,v=1\\ mu \equiv nv \pmod q\\mu \ne nv }}^{q^2}\frac{1}{uv}.
$$
Furthermore, it is shown in~\cite{Zhang}
that 
$$
M_1 = \alpha_q q N + O\(\varphi(q) 2^{\omega(q)} (\log q)^2\).
$$
Thus, it remains to show that
 \begin{equation}
 \label{eq:M2}
M_2 \le N^2 q^{o(1)} .
\end{equation}

Let 
$$
J = \fl{2 \log q}.
$$
Then, changing the order of summation, we obtain
\begin{eqnarray*}
M_2 & = & \varphi(q)\sum_{ u,v=1}^{q^2}\frac{1}{uv} 
 \sum_{\substack{m, n=1\\ mu \equiv nv \pmod q\\mu \ne nv }}^{N} 1\\
 & \le & \varphi(q)\sum_{i,j =0}^J  \sum_{e^i \le u < e^{i+1} } \frac{1}{u} 
  \sum_{e^j \le v < e^{j+1} } \frac{1}{v} 
 \sum_{\substack{m, n=1\\ mu \equiv nv \pmod q\\mu \ne nv }}^{N} 1\\
 & \le & 2 \varphi(q) \sum_{0 \le i \le j \le J}  \sum_{e^i \le u < e^{i+1} } \frac{1}{u} 
  \sum_{e^j \le v < e^{j+1} } \frac{1}{v} 
 \sum_{\substack{m, n=1\\ mu \equiv nv \pmod q\\mu \ne nv }}^{N} 1\\
 & \le & 2 \varphi(q) \sum_{0 \le i \le j \le J} e^{-i-j} \sum_{e^i \le u < e^{i+1} }    \sum_{e^j \le v < e^{j+1} } 
 \sum_{\substack{m, n=1\\ mu \equiv nv \pmod q\\mu \ne nv }}^{N} 1.
 \end{eqnarray*}
 Therefore
 \begin{equation}
\label{eq:M2 and T}
M_2 \le   2 \varphi(q) \sum_{0 \le i \le j \le J} e^{-i - j}  T_{i,j},
\end{equation}
where $T_{i,j}$ is the number of solutions $(m,n,u,v)$ to the congruence
 $$
mu \equiv nv \pmod q, \qquad 1 \le m,n \le N, \ e^i \le u < e^{i+1}, \
e^j \le v < e^{j+1},
$$
with $mu \ne nv$. 

If for a solution $(m,n,u,v)$ we 
write $mu = nv + kq$ with an integer $k$ then  we see that 
$$
1 \le |k| \le  q^{-1} \max\{mu,  nv\}  
 \le  q^{-1} N \max\{e^{i+1},  e^{j+1}\}    =  e^{j+1} N/q.
$$
Thus, there are $O(e^{j}N/q)$ possible values for $k$. 
Clearly  there are at most 
$e^{i+1}$ possible values for  $u$ 
and $N$ possible values $m$. 
Thus the product $nv = mu - kq$ 
can take at most $e^{i+j+2}N^2/q$ possible values and they are 
all of the size $O(Nq^2) = O(q^3)$. Therefore, we see from the bound
on the divisor function~\eqref{eq:phi tau} that
when $m$, $u$ and $k$ are fixed
then $n$ and $v$ can take at most $q^{o(1)}$
possible values.  
Hence
$$
T_{i,j} \le  e^{i+j}N^2 q^{-1 + o(1)}
$$
which after substitution in~\eqref{eq:M2 and T}
gives 
$$
M_2 \le J^2 \varphi(q) N^2 q^{-1 + o(1)}
$$
and the bound~\eqref{eq:M2} follows.  
\end{proof}

 \section{Final Remarks}

As we have mentioned our result is nontrivial for 
$q^{\varepsilon} \le N \le q^{1 -\varepsilon}$. However, 
the author sees no reason why an appropriate asymptotic 
formula cannot hold for even   larger values of $N$, 
say up to $q/2$. It would be interesting to clarify this issue.

   \end{document}